\documentclass[11pt]{article}
\usepackage{amsmath}
\usepackage{amsthm}
 \usepackage{amssymb}
\usepackage{leftidx}
\usepackage{graphicx}
\usepackage{amsfonts}
\usepackage{graphicx, color, epstopdf}
\usepackage{cite}
\usepackage{enumerate}
\usepackage{enumitem}

\newtheorem{theorem}{Theorem}[section]
\newtheorem{lemma}{Lemma}[section]
\newtheorem{proposition}{Proposition}[section]
\newtheorem{corollary}{Corollary}[section]

\newtheorem{remark}{Remark}[section]
\newtheorem{definition}{Definition}[section]

\newcommand{\bey}{\begin{eqnarray}}
\newcommand{\eey}{\end{eqnarray}}
\newcommand{\beyy}{\begin{eqnarray*}}
\newcommand{\eeyy}{\end{eqnarray*}}
\newcommand{\nn}{\nonumber}
\newcommand{\beq}{\begin{equation}}
\newcommand{\eeq}{\end{equation}}
\newcommand{\beqyy}{\begin{equation*}}
\newcommand{\eeqyy}{\end{equation*}}

\newcommand{\fd}[1]{\partial^{#1}_t}
\newcommand{\dfd}[1]{\partial^{#1}_\tau}

\newcommand{\diff}{\triangledown_{\!\tau}}
\newcommand{\defeq}{:=}
\newcommand{\zd}{\,\mathrm{d}}
\newcommand{\abs}[1]{\left|#1\right|}

\newcommand{\bra}[1]{\left(#1\right)}
\newcommand{\brab}[1]{\big(#1\big)}
\newcommand{\braB}[1]{\Big(#1\Big)}

\newcommand{\kbra}[1]{\left[#1\right]}

\newcommand{\mynormb}[1]{\big\|#1\big\|}

\title{Positive definiteness of real quadratic forms
resulting from the variable-step approximation of convolution operators}
\author{Hong-lin Liao\thanks{Department of Mathematics,
Nanjing University of Aeronautics and Astronautics,
Nanjing 211106, P. R. China. E-mails: {\tt liaohl@nuaa.edu.cn}.
The work of this author is supported by a grant 1008-56SYAH18037 from NUAA Scientific Research Starting Fund of Introduced Talent.}
\quad Tao Tang\thanks{
    Division of Science and Technology, BNU-HKBU United International College,
    Zhuhai, Guangdong Province, China \& SUSTech International Center for Mathematics, Shenzhen, China
    Email: {\tt ttang@uic.edu.cn}. The research of this author is supported by science challenge project (TZ2018001) and NSF of China under grant numbers 11731006 and K20911001.}
\quad Tao Zhou\thanks{NCMIS \& LSEC, Institute of Computational Mathematics and Scientific/Engineering Computing,
Academy of Mathematics and Systems Science, Chinese Academy of Sciences, Beijing, 100190,
P. R. China. Email: {\tt tzhou@lsec.cc.ac.cn}. This work of this author is partially supported by NSFC (11822111 and 11688101), Science Challenge Project (TZ2018001),  National Key Basic Research Program (No. 2018YFB0704304), and Youth Innovation promotion Association of CAS.}}
\date{\today}

\begin{document}

\maketitle

\begin{abstract}
The positive definiteness of real quadratic forms
with convolution structures plays an important role in stability analysis for time-stepping schemes for nonlocal operators.
In this work, we present a novel analysis tool to handle discrete convolution kernels resulting from
variable-step approximations for convolution operators. More precisely, for
a class of discrete convolution kernels relevant to variable-step time discretizations,
we show that the associated quadratic form is positive definite under some easy-to-check algebraic conditions. Our proof is based on an elementary constructing strategy using the properties of
discrete orthogonal convolution kernels and complementary convolution kernels. To the best of our knowledge, this is the first general result on simple algebraic conditions for the positive definiteness of variable-step discrete convolution kernels. Using the unified theory, the stability for some simple non-uniform time-stepping schemes can be obtained in a straightforward way.
\end{abstract}


\newpage

\section{Introduction}
\setcounter{equation}{0}

In stability and convergence analysis of time-stepping schemes for integro-differential equations or
time-fractional partial differential equations (see e.g.\cite{ChenXuCaoZhou:2018,Fairweather:1993,LinXu:2007,
LiXu:2013,LubichSloanThomee:1996,PaniFairweather:2002}),
we usually resort to the positive (semi-)definite property of the real quadratic form with convolution structure
\beq
\sum_{k=1}^nw_k\sum_{j=1}^ka_{k-j}w_j\quad\text{for any sequence $\{w_{1}, w_{2}, \cdots, w_{n}\}$,}
\eeq
where the real sequence $\{a_{0},a_{1},\cdots,a_{n},\cdots\}$ is generated by a
time approximation with the uniform time-step of certain convolution integrals such as the Riemann-Liouville fractional integrals
and the Caputo fractional derivatives.

For \textit{uniform} time-stepping approximations, the semi-positive definiteness of the above real quadratic form can be verified by
applying a classical result due to Toeplitz and Carath\'{e}odory \cite[p.18]{GrenanderSzego:1958}. More precisely, if $\{a_{0},a_{1},\cdots,a_{n},\cdots\}$  is a sequence of real numbers such that
\beq \label{eq1a}
\hat{a}(z)=\sum_{k=0}^{\infty}a_nz^k\;\;\text{is analytic in the open unit disk $\mathcal{D}_z:=\{z\in \mathbb{C}: \abs{z}<1\}$,}
\eeq
then the corresponding real quadratic form is positive semi-definite if and only if
\beq \label{eq1b}
\mathrm{Re}\kbra{\hat{a}(z)}\ge0\qquad\text{for any $z\in \mathcal{D}_z$.}
\eeq
This result was widely used in numerical analysis of integro-differential equations, see e.g. \cite{LiXu:2013,Lubich_I, Lubich_II, LubichSloanThomee:1996,PaniFairweather:2002}.

In analyzing finite difference schemes for nonlinear integro-differential
equations, L\'{o}pez-Marcos  \cite[Proposition 5.2]{LopezMarcos:1990} presented a corollary of the above result with
the following sufficient conditions for the real sequence $\{a_{0},a_{1},\cdots,a_{n},\cdots\}:$
\begin{align}\label{cond: LopezMarcos 1990}
a_j\ge0,\quad a_{j-1}\ge a_{j},\quad a_{j-1}-a_{j}\ge a_{j}-a_{j+1}.
\end{align}
Notice that the above criterions are rather simple and easy-to-check. Actually, it has been shown to be powerful in stability and convergence
analysis of \textit{uniform} time-stepping approximations for integro-differential and time-fractional differential problems, see e.g., \cite{CuestaPalencia:2003,JiLiaoZhang:2019,SunWu:2006,Tang:1993}.

However, in most practical applications, one may need to use nonuniform/variable time-stepping schemes. This is the case when one wants to capture the multi-scale behaviors in time, and/or one needs to handle solution singularities in time. In these situations, it is natural to investigate (semi-)positive definiteness of the following real quadratic form with convolution structure:
\begin{align}\label{eq: real quadratic form}
\sum_{k=1}^nw_k\sum_{j=1}^ka^{(k)}_{k-j}w_j\quad\text{for any sequence $\{w_{1}, w_{2}, \cdots, w_{n}\}$,}
\end{align}
where the associated (variable) discrete convolution kernels $\big\{a^{(n)}_{n-k}\big\}_{k=1}^n$ may be resulted from
variable time-stepping approximations of the Caputo fractional derivative and Riemann--Liouville (convolution) integrals, see e.g., \cite{JiLiaoGongZhang:2019,JiLiaoZhang:2019,LiaoLiZhang:2018,LiaoMcLeanZhang:2019,LiaoTangZhou:2019,
McLeanMustapha:2007,Mustapha:2011, McLeanThomeeWahlbin:1996}. Here the superscript $(n)$
represents the discrete time-level $t_n$ and the value of $a^{(n)}_{j}$ may vary at different time levels $t_n.$

It is noted that the extension of (\ref{eq1a})-(\ref{eq1b}) or (\ref{cond: LopezMarcos 1990}) to nonuniform time-stepping schemes is by no means trivial. In other words, one can not apply the criterions in \cite{GrenanderSzego:1958, LopezMarcos:1990} to examine
the positive definiteness  of the quadratic form \eqref{eq: real quadratic form}. We also remark that one may consider the so-called positive-semidefinite-preserving approach \cite{JiLiaoGongZhang:2019,McLeanThomeeWahlbin:1996,McLeanMustapha:2007,Mustapha:2011} to deal with the positive definiteness of  \eqref{eq: real quadratic form} as the associated discrete convolution kernels are discrete analogy of the positive definite continuous kernels. However, these specially designed discrete kernels are of very special form and the corresponding theory lacks application generality.

The aim of this work is to present some sufficient and easy-to-check criterions
on a wide class of discrete convolution kernels $\big\{a^{(n)}_{n-k}\big\}_{k=1}^n$
to ensure the semi-positive definiteness of \eqref{eq: real quadratic form}. For notation simplicity, we shall simply use $\mathcal{B}^{(n)}_{n-k}$
to represent a set of discrete convolution kernels
$\big\{\mathcal{B}^{(n)}_{n-k}\big\}_{k=1}^n$ or
$\big\{\mathcal{B}^{(n)}_{0},\mathcal{B}^{(n)}_{1},\cdots, \mathcal{B}^{(n)}_{n-1}\big\}.$
We say that the discrete convolution kernels $\mathcal{B}^{(n)}_{n-k}$
are positive (semi-)definite if the associated real quadratic form is positive (semi-)definite.

We now state our main theorem as follows.

\begin{theorem}\label{thm: positive definite quadratic form}
For fixed $n\ge2$, the discrete convolution kernels $a^{(n)}_{n-k}$ are positive definite if the following conditions are satisfied
\beyy
{\bf C1.} \qquad &a^{(n)}_{j}>0 & \mbox{for} \;\; 0\le j \le n-1; \\
{\bf C2.} \qquad &a^{(n-1)}_{j-1}>a^{(n)}_{j} & \mbox{for} \;\; 1\le j \le n-1;\\
{\bf C3.} \qquad &a^{(n-1)}_{j-1}a^{(n)}_{j+1}\ge a^{(n-1)}_{j}a^{(n)}_{j}\quad & \mbox{for} \;\; 1\le j \le n-2;\\
{\bf C4.} \qquad &a^{(n)}_{j-1}\ge a^{(n)}_{j} &\mbox{for} \;\; 1\le j \le n-1.
\eeyy
\end{theorem}
\vskip .25cm

To the best of our knowledge, Theorem \ref{thm: positive definite quadratic form}  is the first result with simple algebraic conditions for the positive definiteness of variable convolution coefficients. Moreover, the conditions stated in Theorem \ref{thm: positive definite quadratic form} do not have explicit connection with the associated continuous kernel.

The proof of  Theorem \ref{thm: positive definite quadratic form} will be presented in Section \ref{sect3}.
Our proof is based on two novel discrete tools. The first one is the so-called \emph{discrete orthogonal convolution} (DOC) kernels $\theta_{n-k}^{(n)}$,
which satisfy the following \emph{discrete orthogonal identity}
\beq\label{eq: orthogonal identity0}
\sum_{j=k}^{n}\theta_{n-j}^{(n)}a^{(j)}_{j-k}\equiv\delta_{nk}
\quad\text{for any $\;1\leq k\le n$,}
\eeq
where $\delta_{nk}$ is the Kronecker delta symbol. We will demonstrate in Lemma \ref{lem: Mutual orthogonality} that positive (semi-)definiteness of the DOC kernels $\theta_{n-k}^{(n)}$ implies the positive (semi-)definiteness of $a_{n-k}^{(n)}.$ Then it remains to show the positive definiteness of the DOC kernels $\theta_{n-k}^{(n)}$
under conditions \textbf{C1}-\textbf{C4}. The proof will rely on our second discrete tool, i.e.,
the so-called \emph{discrete complementary convolution} (DCC) kernels $p_{n-k}^{(n)}$ which satisfy
the following \emph{discrete complementary identity:}
\beq\label{eq: complementary identity0}
\sum_{j=k}^{n}p_{n-j}^{(n)}a^{(j)}_{j-k}\equiv1\quad\text{for any $\;1\leq k\le n$.}
\eeq

The rest of the paper is organized in the following way. In Section \ref{sect2}, we present some useful theory for DOC and DCC.
Then we will provide the rigorous proof for Theorem \ref{thm: positive definite quadratic form} in Section \ref{sect3}. In Section \ref{sect4},
we present some applications of Theorem \ref{thm: positive definite quadratic form}
to variable time-stepping approximations for the Caputo fractional derivative and general convolution integrals. Finally, we provide some concluding remarks in Section 5. In particular, we will discuss possible ways to use the DOC and DCC kernels for handling more general discrete convolution kernels.

\section{On discrete orthogonal and discrete complementary convolution kernels} \label{sect2}
\setcounter{equation}{0}

For generality, we will not restrict our discussions to the specific discrete kernels
satisfying \textbf{C1}-\textbf{C4} as stated in Theorem \ref{thm: positive definite quadratic form}.
Without loss of generality, we shall assume that the discrete convolution kernels satisfy
$a_{j}^{(n)}\neq0$ for $0\le j\le n-1.$
\beq\label{assume: discrete kernels a}
a_{j}^{(n)}=0\;\;\text{for $j\ge n_0$,}\;\;
\text{when $a_{n_0}^{(n)}=0$ appears for some fixed index $2\le n_0\le n-1$.}
\eeq
Notice that this assumption is reasonable in approximating nonlocal operators.

\subsection{Discrete orthogonal convolution kernels}

The discrete orthogonal convolution (DOC) kernels $\theta_{n-j}^{(n)}$ induced by $a_{n-j}^{(n)}$ are defined via a recursive procedure
\beq\label{eq: orthogonal recursive procedure}
\theta_{0}^{(n)}:=\frac{1}{a^{(n)}_{0}}\quad\text{and}\quad
\theta_{n-k}^{(n)}:=-\frac1{a^{(k)}_{0}}\sum_{j=k+1}^{n}\theta_{n-j}^{(n)}a^{(j)}_{j-k}
\quad\text{for $k=n-1, n-2, \cdots, 1$.}
\eeq
It can be verified that the following discrete orthogonal identity holds
\beq\label{eq: orthogonal identity}
\sum_{j=k}^{n}\theta_{n-j}^{(n)}a^{(j)}_{j-k}= \delta_{nk}
\quad\text{for $1\leq k\le n$},
\eeq
where $\delta_{nk}$ is the Kronecker delta symbol.

Note that, under the assumption \eqref{assume: discrete kernels a},
the DOC kernels $\theta_{n-j}^{(n)}$ are uniquely determined by the original discrete convolution kernels $a^{(n)}_{n-k}$.
This type of DOC kernels was originally proposed in \cite{LiaoZhang:2019} for analyzing the $L^2$-stability of
variable-step BDF2 schemes.

Below we will present a general theory for the DOC kernels (\ref{eq: orthogonal recursive procedure}), which establishes the positive definiteness equivalence between $\theta_{n-j}^{(n)}$ and its associated kernel $a^{(n)}_{n-k}$.

\begin{lemma}\label{lem: Mutual orthogonality}
The DOC kernel $\theta_{n-j}^{(n)}$ defined by \eqref{eq: orthogonal recursive procedure}
and its associated kernel $a^{(n)}_{n-k}$ are mutually orthogonal, i.e.,
\beq\label{eq: Mutual orthogonal identity}
\sum_{j=k}^na^{(n)}_{n-j}\theta_{j-k}^{(j)}=\delta_{nk}, \qquad  \sum_{j=k}^{n}\theta_{n-j}^{(n)}a^{(j)}_{j-k}=\delta_{nk},
\quad \text{for any $1\le k\le n$.}
\eeq
Consequently, $\theta_{n-j}^{(n)}$ is positive (semi-)definite if and only if  $a^{(n)}_{n-k}$ is positive (semi-)definite.
\end{lemma}

\begin{proof}
For any fixed index $n\ge1$ and any sequence $\{w_{k}\}_{k=1}^n$, a sequence $\{V_{k}\}_{k=1}^n$ is defined by
\begin{align}\label{def: V-w}
V_j:=\sum_{k=1}^{j}a_{j-k}^{(j)}w_k, \quad 1\le j\le n.
\end{align}
Consider the matrix $\mathrm{A}=\brab{\mathrm{a}_{ij}}_{n\times n}$ whose elements are given by
\begin{equation}
\mathrm{a}_{ij}:=a_{i-j}^{(i)} \,\,\, \textmd{for} \,\,\, i\le j  \quad
\textmd{and} \quad  \mathrm{a}_{ij}:=0 \,\,\, \mbox{for} \,\,\, i> j.
\end{equation}
Then, under the general setting \eqref{assume: discrete kernels a}, the non-singularity of $\mathrm{A}$ implies that
$\{V_j\}_{j=1}^n$ is non-zero if and only if $\{w_j\}_{j=1}^n$ is non-zero.

Now multiplying both sides of \eqref{def: V-w} by $\theta_{n-j}^{(n)}$ and summing $j$ up from $j=1$ to $j=n$ lead to
\beyy
\sum_{j=1}^n\theta_{n-j}^{(n)}V_j &=&
\sum_{j=1}^n\theta_{n-j}^{(n)}\sum_{k=1}^{j}a_{j-k}^{(j)}w_k= \sum_{k=1}^{n}w_k
\sum_{j=k}^n\theta_{n-j}^{(n)}a_{j-k}^{(j)} \nn \\
&=& \sum_{k=1}^{n}w_k\delta_{nk}=w_n,
\eeyy
where the discrete orthogonal identity \eqref{eq: orthogonal identity} is used in the above derivation.
Thus we obtain a formula for $w_j$ via $V_k$:
\begin{align}\label{eq: w-V orthogonal relationship}
w_j=\sum_{k=1}^j\theta_{j-k}^{(j)}V_k, \qquad  1\le j\le n,
\end{align}
which can be regarded as an inverse for (\ref{def: V-w}).
Using the definition (\ref{def: V-w}), and multiplying both sides of (\ref{eq: w-V orthogonal relationship}) by
$a^{(n)}_{n-j}$ and summing $j$ up from $j=1$ to $j=n$, give
\beqyy 
V_n = \sum_{j=1}^na^{(n)}_{n-j}w_j=
\sum_{j=1}^na^{(n)}_{n-j}\sum_{k=1}^j\theta_{j-k}^{(j)}V_k = \sum_{k=1}^nV_k\sum_{j=k}^na^{(n)}_{n-j}\theta_{j-k}^{(j)}, \qquad n\ge1.
\eeqyy
Using the arbitrariness of $\{w_k\}$ (which in turn yields  arbitrariness of $\{V_k\}$ due to their equivalence) yields
\begin{align*}
\sum_{j=k}^na^{(n)}_{n-j}\theta_{j-k}^{(j)}= \delta_{nk},
\quad 1\le k\le n.
\end{align*}
This yields the desired result \eqref{eq: Mutual orthogonal identity}. Furthermore,
the positive (semi-)definiteness of the DOC kernels $\theta_{n-k}^{(n)}$
(or the original kernels $a^{(n)}_{n-j}$) directly follows from the following identity
\beq \label{2eqc0}
\sum_{k=1}^nw_k\sum_{j=1}^ka^{(k)}_{k-j}w_j=\sum_{k=1}^nV_k\sum_{j=1}^k\theta^{(k)}_{k-j}V_j,
\eeq
which can be obtained by using \eqref{def: V-w} and \eqref{eq: w-V orthogonal relationship}.
\end{proof}

Lemma \ref{lem: Mutual orthogonality} implies that the positive definiteness for
the discrete convolution kernel $a_{n-k}^{(n)}$ is equivalent to that for its DOC kernel $\theta_{n-j}^{(n)}$.

Notice that by Definition \ref{eq: orthogonal recursive procedure}, one can apply
the conditions \textbf{C1} and \textbf{C3} to show that
\beyy
&& \theta_{0}^{(n)}=1/{a^{(n)}_{0}}>0, \\
&&  \theta_{1}^{(n)}=\,-\frac{a^{(n)}_{1}}{a^{(n-1)}_{0}}\theta_{0}^{(n)}<0,\\
&&\theta_{2}^{(n)}=\,-\frac1{a^{(n-2)}_{0}}\bra{\theta_{0}^{(n)}a^{(n)}_{2}+\theta_{1}^{(n)}a^{(n-1)}_{1}} =-\frac{a^{(n-1)}_{1}}{a^{(n-2)}_{0}}\theta_{0}^{(n)}\braB{\frac{a^{(n)}_{2}}{a^{(n-1)}_{1}}-\frac{a^{(n)}_{1}}{a^{(n-1)}_{0}}}<0.
\eeyy

 Below we will introduce two auxiliary sequences $\psi_j^{(m)}$ and $\chi_{\ell}^{(k)}$. The purpose of doing that is to derive an explicit formula for $\theta_{n-j}^{(n)}$, i.e., not via the recurrence generation.
 .
\begin{definition}\label{def: auxiliary sequence Psi}
For any fixed $m\ge2$, we define the auxiliary sequence $\psi_{j}^{(m)}$ by
\bey
\psi^{(m)}_0 :=1/{a^{(n)}_{0}} \quad\text{and}\quad
\psi_{j}^{(m)}:={a_{j}^{(m)}}/{a_{j-1}^{(m-1)}}, \;\quad\; 1\le j\le m-1, \label{2eqc2}
\eey
where as conventional understanding  we set $\psi_{j}^{(m)}=0$ if the denominator $a_{j-1}^{(m-1)}=0$ for $1\le j\le m-1$.
\end{definition}

\begin{definition}\label{def: auxiliary sequence Chi}
For $n\ge3$, we define the auxiliary sequence $\{\chi_{\ell}^{(k)}\,|\, 2\leq \ell \leq n-k-1\}$
for $0\leq k \leq n-3$ using $\psi_{j}^{(m)}$ given by \eqref{2eqc2}:
\begin{itemize}
  \item For $k=0$, we set $\chi_\ell^{(0)}:=1$ \,\,for\,\, $2\leq \ell \leq  n-1$;
  \item For $k=1$, we set
  \beqyy 
  \chi_\ell^{(1)}:=\frac{\psi_{\ell+1}^{(n)}-\psi_{1}^{(n)}}{\psi_{2}^{(n)}-\psi_{1}^{(n)}}
  =\frac{\chi_{\ell+1}^{(0)}\psi_{\ell+1}^{(n)}-\psi_{1}^{(n)}}{\chi_2^{(0)}\psi_{2}^{(n)}-\psi_{1}^{(n)}},\qquad2\leq \ell \leq n-2;
  \eeqyy
  \item For $2\leq k \leq n-3,$ we set
  \beq \label{2eqc4}
  \chi_\ell^{(k)}:=\frac{\chi_{\ell+1}^{(k-1)}\psi_{\ell+1}^{(n-k+1)}-\psi_{1}^{(n-k+1)}}{\chi_2^{(k-1)}
  \psi_{2}^{(n-k+1)}-\psi_{1}^{(n-k+1)}}, \qquad 2\leq \ell \leq n-k-1.
  \eeq
\end{itemize}
As usual, we set $\chi_\ell^{(k)}:=1$ whenever the denominator in (\ref{2eqc4}) is 0.
\end{definition}

We point out that although $\chi_{\ell}^{(k)}$ is defined recursively via the auxiliary sequence $\psi_{j}^{(m)}$ it can be verified that it is uniquely determined by the original discrete kernels $a_{n-k}^{(n)}$,

Now let us deduce the first few entries of DOC kernels by using the auxiliary sequence $\psi_{j}^{(m)}$  and $\chi_{\ell}^{(k)}$. It can be verified tat
\bey
\theta_{1}^{(n)} &=&-\psi_{0}^{(n)}\psi_{1}^{(n)}, \qquad n\ge2, \label{2eqc6}\\
\theta_{2}^{(n)} &=&-\psi_{0}^{(n)}\psi_{1}^{(n-1)}\bra{\psi_{2}^{(n)}-\psi_{1}^{(n)}}, \qquad n\ge3, \label{2eqc7}\\
\theta_{3}^{(n)}&=&-\frac1{a^{(n-3)}_{0}}\bra{\theta_{0}^{(n)}a^{(n)}_{3}+\theta_{1}^{(n)}a^{(n-1)}_{2}
+\theta_{2}^{(n)}a^{(n-2)}_{1}}\nn \\
&=&-\psi_{1}^{(n-2)}\bra{\theta_{0}^{(n)}\psi_{3}^{(n)}\psi_{2}^{(n-1)}
+\theta_{1}^{(n)}\psi_{2}^{(n-1)}+\theta_{2}^{(n)}}\nn \\
&=&-\psi_{0}^{(n)}\psi_{1}^{(n-2)} \bra{\chi_2^{(1)}\psi_{2}^{(n-1)}-\psi_{1}^{(n-1)}}\bra{\chi_2^{(0)}\psi_{2}^{(n)}-\psi_{1}^{(n)}}, \quad n\ge4. \label{2eqc8}
\eey
It is obvious that above gives an explicit formulas for the DOC kernels $\theta_{j}^{(n)}$. We now extend the above observations to the following lemma.

\begin{lemma}\label{lem: DOC explicit formula}
For fixed $n\ge2$, the DOC kernels $\theta_{n-j}^{(n)}$ defined in \eqref{eq: orthogonal recursive procedure} satisfy
\beq\label{eq: orthogonal kernel formula}
\theta_{j}^{(n)}=-\psi_{0}^{(n)}\psi^{(n-j+1)}_1\prod_{\ell=n-j+2}^{n}\brab{\chi^{(n-\ell)}_{2}\psi_{2}^{(\ell)}-\psi_{1}^{(\ell)}},
\quad\; 1\leq j\leq n-1,
\eeq
where $\psi_{j}^{(m)}$  and $\chi_{\ell}^{(k)}$ are defined by Definitions \ref{def: auxiliary sequence Psi} and \ref{def: auxiliary sequence Chi}.
\end{lemma}
\begin{proof}
It follows from (\ref{2eqc6})-(\ref{2eqc8}) that \eqref{eq: orthogonal kernel formula} holds for $j=1,2,3$.
For notation simplicity, in what follows we set
\[
\Psi^{(m)}=\chi^{(n-m)}_{2}\psi_{2}^{(m)}-\psi_{1}^{(m)}.
\]
For $n\ge 4$, it follows from Definition \ref{def: auxiliary sequence Chi} that
\begin{align}\label{eq: sequence PSI}
\chi_{\ell+1}^{(n-m)}\psi_{\ell+1}^{(m)}-\psi_{1}^{(m)}=\chi_\ell^{(n-m+1)}\Psi^{(m)}\quad\text{for $\ell=2,\cdots, n-2$.}
\end{align}
We will perform the proof by induction. To begin, we assume that the formula \eqref{eq: orthogonal kernel formula} is true for
$1\le j \le k-1$, i.e.,
\begin{align}\label{eq: induction hypothesis}
\theta_{j}^{(n)}=\psi_0^{(n)}\psi^{(n-j+1)}_1\prod_{\ell=n-j+2}^{n}\Psi^{(\ell)},\quad
\quad\text{for $1\leq j\leq k-1$,}
\end{align}
and we will prove \eqref{eq: orthogonal kernel formula} is satisfied for $j=k$. It follows from the definition \eqref{2eqc2} that
\beqyy
a_{j-k}^{(j)}=\psi_{j-k}^{(j)}\psi_{j-k-1}^{(j-1)}\cdots\psi_{1}^{(k+1)}a_{0}^{(k)}
=a_{0}^{(k)}\prod_{\ell=k+1}^{j}\psi_{\ell-k}^{(\ell)}, \qquad\text{for $1\le k\le j$.}
\eeqyy
Replacing the index $k$ by $n-k$ gives
\beqyy
a_{j-n+k}^{(j)}=\psi_{j-n+k}^{(j)}\psi_{j-n+k-1}^{(j-1)}\cdots\psi_{1}^{(n-k+1)}a_{0}^{(n-k)}
=a_{0}^{(n-k)}\prod_{\ell=n-k+1}^{j}\psi_{\ell-n+k}^{(\ell)}\,.
\eeqyy
It follows from Definition \eqref{eq: orthogonal recursive procedure} and the induction hypothesis \eqref{eq: induction hypothesis} that
\beyy
\theta_{k}^{(n)}&=& -\frac1{a^{(n-k)}_{0}}\psi_{0}^{(n)}a^{(n)}_{k}
-\frac1{a^{(n-k)}_{0}}\sum_{j=n-k+1}^{n-1}\theta_{n-j}^{(n)}a^{(j)}_{j-n+k}
\nn \\
&=&-\psi_{0}^{(n)}\prod_{\ell=n-k+1}^{n}\psi_{\ell-n+k}^{(\ell)}
-\sum_{j=n-k+1}^{n-1}\theta_{n-j}^{(n)}\prod_{\ell=n-k+1}^{j}\psi_{\ell-n+k}^{(\ell)}\nn \\
&=&-\psi_{0}^{(n)}\prod_{\ell=n-k+1}^{n}\psi_{\ell-n+k}^{(\ell)}+\theta_0^{(n)}\psi^{(n)}_1\prod_{\ell=n-k+1}^{n-1}\psi_{\ell-n+k}^{(\ell)}
-\sum_{j=n-k+1}^{n-2}\theta_{n-j}^{(n)}\prod_{\ell=n-k+1}^{j}\psi_{\ell-n+k}^{(\ell)}\nn \\
&=&-\psi_{0}^{(n)}\brab{\psi_{k}^{(n)}-\psi^{(n)}_{1}}\prod_{\ell=n-k+1}^{n-1}\psi_{\ell-n+k}^{(\ell)}
-\sum_{j=n-k+1}^{n-2}\theta_{n-j}^{(n)}\prod_{\ell=n-k+1}^{j}\psi_{\ell-n+k}^{(\ell)}\nn \\
&=& -\psi_{0}^{(n)}\chi_{k-1}^{(1)}\Psi^{(n)}\prod_{\ell=n-k+1}^{n-1}\psi_{\ell-n+k}^{(\ell)}
-\sum_{j=n-k+1}^{n-2}\theta_{n-j}^{(n)}\prod_{\ell=n-k+1}^{j}\psi_{\ell-n+k}^{(\ell)}\,, 
\eeyy
where in the past step the equality \eqref{eq: sequence PSI} is used with $m=n$ and $\ell=k-1$.
Furthermore, with the help of the hypothesis \eqref{eq: induction hypothesis}, we have
\beyy
\theta_{k}^{(n)}&=&-\bra{\psi_{0}^{(n)}\psi_{k-1}^{(n-1)}\chi_{k-1}^{(1)}\Psi^{(n)}
+\theta_{2}^{(n)}}\prod_{\ell=n-k+1}^{n-2}\psi_{\ell-n+k}^{(\ell)}
-\sum_{j=n-k+1}^{n-3}\theta_{n-j}^{(n)}\prod_{\ell=n-k+1}^{j}\psi_{\ell-n+k}^{(\ell)}\nn \\
&=&-\psi_{0}^{(n)}\bra{\chi_{k-1}^{(1)}\psi_{k-1}^{(n-1)}-\psi^{(n-1)}_1}\Psi^{(n)}\prod_{\ell=n-k+1}^{n-2}\psi_{\ell-n+k}^{(\ell)}
-\sum_{j=n-k+1}^{n-3}\theta_{n-j}^{(n)}\prod_{\ell=n-k+1}^{j}\psi_{\ell-n+k}^{(\ell)}\nn \\
&=&-\psi_{0}^{(n)}\chi_{k-2}^{(2)}\Psi^{(n-1)}\Psi^{(n)}\prod_{\ell=n-k+1}^{n-2}\psi_{\ell-n+k}^{(\ell)}
-\sum_{j=n-k+1}^{n-3}\theta_{n-j}^{(n)}\prod_{\ell=n-k+1}^{j}\psi_{\ell-n+k}^{(\ell)}\,, 
\eeyy
where in the last step \eqref{eq: sequence PSI} is used with $m=n-1$ and $\ell=k-2$.
Repeating the above process gives
\beyy
\theta_{k}^{(n)}
&=&-\psi_{0}^{(n)}\chi_{k-2}^{(2)}\prod_{\ell=n-1}^{n}\Psi^{(\ell)}\prod_{\ell=n-k+1}^{n-3}\psi_{\ell-n+k}^{(\ell)}
-\sum_{j=n-k+1}^{n-4}\theta_{n-j}^{(n)}\prod_{\ell=n-k+1}^{j}\psi_{\ell-n+k}^{(\ell)}\nn \\
&=&-\psi_{0}^{(n)}\chi_{k-3}^{(3)}\prod_{\ell=n-2}^{n}\Psi^{(\ell)}\prod_{\ell=n-k+1}^{n-4}\psi_{\ell-n+k}^{(\ell)}
-\sum_{j=n-k+1}^{n-5}\theta_{n-j}^{(n)}\prod_{\ell=n-k+1}^{j}\psi_{\ell-n+k}^{(\ell)}\nn \\
&=&\,\cdots\cdots \nn \\
&=&-\psi_{0}^{(n)}\chi_{2}^{(k-2)}\prod_{\ell=n-k+3}^{n}\Psi^{(\ell)}\prod_{\ell=n-k+1}^{n-k+2}\psi_{\ell-n+k}^{(\ell)}
-\theta_{k-1}^{(n)}\psi_{1}^{(n-k+1)}\nn \\
&=&-\psi_{0}^{(n)}\psi_{1}^{(n-k+1)}\brab{\chi_{2}^{(k-2)}\psi^{(n-k+2)}_2-\psi^{(n-k+2)}_1}\prod_{\ell=n-k+3}^{n}\Psi^{(\ell)}\nn \\
&=&-\psi_{0}^{(n)}\psi_{1}^{(n-k+1)}\prod_{\ell=n-k+2}^{n}\Psi^{(\ell)}\,. 
\eeyy
This proves \eqref{eq: orthogonal kernel formula} for $j=k$. The mathematical induction completes the proof for
\eqref{eq: orthogonal kernel formula}.
\end{proof}

The next lemma will examine the signs of the DOC kernels $\theta_{n-k}^{(n)}.$

\begin{lemma}\label{lem: sign of DOC kernels}
Assume that the discrete convolution kernels $a_{n-k}^{(n)}$ satisfy the conditions \textbf{C1}-\textbf{C3}. Then
for any $n\ge2$, the DOC kernels $\theta_{n-k}^{(n)}$ in \eqref{eq: orthogonal recursive procedure} satisfy
\begin{itemize}
\item
the sign property:
\bey
&& \theta_{0}^{(n)}>0,\quad \theta_{1}^{(n)}<0\quad \text{and}\quad \theta_{j}^{(n)}\le0\;\quad \text{for $2\le j\le n-1$,}\label{2eqc15} \\
&& \sum_{j=1}^n\theta_{n-j}^{(n)}>0; \label{2eqc16}
\eey
\item
and the convolution quadratic inequality, i.e., for any sequence $\{V_j\}_{j=1}^n$,
\beq \label{2eqc17}
2\sum_{k=1}^nV_k\sum_{j=1}^k\theta_{k-j}^{(k)}V_j
\ge\sum_{k=1}^n\bigg(\sum_{j=1}^k\theta_{k-j}^{(k)}+\sum_{j=k}^n\theta_{j-k}^{(j)}\bigg)V_k^2.
\eeq
\end{itemize}
\end{lemma}

\begin{proof} Using the condition \textbf{C1} gives $\theta_{0}^{(n)}=1/a_{0}^{(n)}>0$. It follows the
conditions \textbf{C1} and \textbf{C3} that
\beqyy
\frac{a^{(m)}_{j+1}}{a^{(m-1)}_{j}}\ge\frac{a^{(m)}_{j}}{a^{(m-1)}_{j-1}}\, , \qquad\text{for $1\le j\le m-2$}.
\eeqyy
Consequently, using Definition \ref{def: auxiliary sequence Psi} gives
\[
\psi_{m-1}^{(m)}\ge\cdots\ge\psi_{2}^{(m)}\ge\psi_{1}^{(m)}>0\quad\text{for $2\le m \le n$.}
\]
It follows from Definition \ref{def: auxiliary sequence Chi} that $\chi_{\ell}^{(0)}=1$ for any $\ell\ge2$, and
\[
\chi_{\ell+1}^{(1)}=\frac{\psi_{\ell+2}^{(n)}-\psi_{1}^{(n)}}{\psi_{2}^{(n)}-\psi_{1}^{(n)}}
\ge\frac{\psi_{\ell+1}^{(n)}-\psi_{1}^{(n)}}{\psi_{2}^{(n)}-\psi_{1}^{(n)}}=\chi_\ell^{(1)}\ge1\quad\text{for $\ell\ge2$.}
\]
By a simple induction for $k=n-1,n-2,\cdots,4$, it is not difficult to check that
\[
\chi_2^{(n-k)}\psi_{2}^{(k)}-\psi_{1}^{(k)}\ge\psi_{2}^{(k)}-\psi_{1}^{(k)}\ge0 .
\]
Consequently,
\beqyy
\chi_{\ell+1}^{(n-k+1)}
=\frac{\chi_{\ell+2}^{(n-k)}\psi_{\ell+2}^{(k)}-\psi_{1}^{(k)}}{\chi_2^{(n-k)}\psi_{2}^{(k)}-\psi_{1}^{(k)}}
\ge\frac{\chi_{\ell+1}^{(n-k)}\psi_{\ell+1}^{(k)}-\psi_{1}^{(k)}}{\chi_2^{(n-k)}\psi_{2}^{(k)}-\psi_{1}^{(k)}}
=\chi_{\ell}^{(n-k+1)}\ge1
\eeqyy
for $\ell=2,\cdots,k-3$. Thus using Lemma \ref{lem: DOC explicit formula} yields that $\theta_{1}^{(n)}<0$ and
\beqyy
\theta_{n-j}^{(n)}\le-\theta_{0}^{(n)}\psi^{(j+1)}_1\prod_{\ell=j+2}^{n}\brab{\psi_{2}^{(\ell)}-\psi_{1}^{(\ell)}}\le 0
\quad\text{for $1\leq j\leq n-2$.}
\eeqyy
This proves (\ref{2eqc15}). To prove (\ref{2eqc16}), we first let
\[
\sigma_j:=\sum_{k=1}^j\theta_{j-k}^{(j)}.
\]
The orthogonal identity in \eqref{eq: Mutual orthogonal identity} implies that
\beq \label{2eqc20}
\sum_{j=1}^ma_{m-j}^{(m)}\sigma_j=\sum_{j=1}^ma_{m-j}^{(m)}\sum_{k=1}^j\theta_{j-k}^{(j)}
=\sum_{k=1}^m\sum_{j=k}^ma_{m-j}^{(m)}\theta_{j-k}^{(j)}\equiv1\quad \text{for any $m\ge1$}.
\eeq
Taking $m=1$ yields $\sigma_1=1/a_{0}^{(1)}>0$. For $m\ge2$, taking $m=\ell$ and $\ell-1$ in (\ref{2eqc20}) gives
\beq\label{eq: sign of DOC summation}
\sigma_\ell\equiv\frac{1}{a_{0}^{(\ell)}}\sum_{j=1}^{\ell-1}\brab{a_{\ell-j-1}^{(\ell-1)}-a_{\ell-j}^{(\ell)}}\sigma_j\quad \text{for any $\ell\ge2$}.
\eeq
The  condition \textbf{C2} ensures the positivity of the coefficient on the right-hand side of (\ref{eq: sign of DOC summation}).
Thus a simple induction shows that $\sigma_\ell>0\quad\text{for $\ell\ge2$}$
by taking $\ell=2,\cdots,n$ in \eqref{eq: sign of DOC summation} successively. This confirms (\ref{2eqc16}).

Now applying (\ref{2eqc15})-(\ref{2eqc16}) together with the conditions \textbf{C1}-\textbf{C3} yield
\begin{align*}
2V_k\sum_{j=1}^k\theta_{k-j}^{(k)}V_j=&\,2\theta_{0}^{(k)}V_k^2+2\sum_{j=1}^{k-1}\theta_{k-j}^{(k)}V_kV_j\\
\ge&\,2\theta_{0}^{(k)}V_k^2+\sum_{j=1}^{k-1}\theta_{k-j}^{(k)}\bra{V_k^2+V_j^2}\\
=&\,V_k^2\sum_{j=1}^k\theta_{k-j}^{(k)}+\sum_{j=1}^k\theta_{k-j}^{(k)}V_j^2\quad \text{for $k\ge1$.}
\end{align*}
By summing up $k$ from $k=1$ to $n$ and exchanging the order of summation, one gets
\begin{align*}
2\sum_{k=1}^nV_k\sum_{j=1}^k\theta_{k-j}^{(k)}V_j
\ge&\,\sum_{k=1}^nV_k^2\sum_{j=1}^k\theta_{k-j}^{(k)}+\sum_{k=1}^n\sum_{j=1}^k\theta_{k-j}^{(k)}V_j^2\\
=&\,\sum_{k=1}^nV_k^2\sum_{j=1}^k\theta_{k-j}^{(k)}+\sum_{k=1}^nV_k^2\sum_{j=k}^n\theta_{j-k}^{(j)}\quad \text{for $n\ge1$,}
\end{align*}
which yields the desired result (\ref{2eqc17}).
\end{proof}

\subsection{The discrete complementary convolution kernels}

The lower bound in (\ref{2eqc17}) motivates us to define a new class of discrete kernels $p_{n-j}^{(n)}$ by using the DOC kernels $\theta_{n-j}^{(n)}$:
\begin{align}\label{eq: summing orthogonal procedure}
p_{n-k}^{(n)}:=\sum_{j=k}^{n}\theta_{j-k}^{(j)}\quad\text{for $1\leq k\le n$.}
\end{align}
Obviously, the uniqueness of $\theta_{n-j}^{(n)}$ guarantees the uniqueness of $p_{n-j}^{(n)}$.
Notice that
\begin{align}\label{eq: orthogonal-complementary relation}
\theta_{0}^{(n)}=p_{0}^{(n)}\quad\text{and}\quad
\theta_{n-k}^{(n)}=p_{n-k}^{(n)}-p_{n-k-1}^{(n-1)}\quad\text{for $1\leq k\le n-1$.}
\end{align}
Inserting the above equations into the discrete orthogonal identity \eqref{eq: orthogonal identity} gives
\begin{align}\label{eq: recursive complementary}
\sum_{j=k}^{n}p_{n-j}^{(n)}a^{(j)}_{j-k}
\equiv\sum_{j=k}^{n-1}p_{n-1-j}^{(n-1)}a^{(j)}_{j-k}+\delta_{nk}\quad\text{for $\;\;1\leq k\le n$.}
\end{align}
By setting $\Xi_{k}^{(n)}:=\sum_{j=k}^{n}p_{n-j}^{(n)}a^{(j)}_{j-k}$, one obtains
\begin{align*}
\Xi_{k}^{(k)}=1\quad\text{and}\quad \Xi_{k}^{(n)}=\Xi_{k}^{(n-1)}\quad\text{for $\;\;1\leq k\le n-1$.}
\end{align*}
Thus a simple induction yields the following discrete complementary identity
\begin{align}\label{eq: new complementary identity}
\sum_{j=k}^{n}p_{n-j}^{(n)}a^{(j)}_{j-k}\equiv1\quad\text{for $\;\;1\leq k\le n$.}
\end{align}
This implies that the discrete kernels $p_{n-j}^{(n)}$ in \eqref{eq: summing orthogonal procedure} are
complementary to the original kernels $a_{n-j}^{(n)}$, which also explains why  $p_{n-j}^{(n)}$ is called the
discrete complementary convolution kernels.

We present in the next lemma an explicit formulation of the DCC kernels $p_{n-k}^{(n)}$,
which only relies on the original discrete convolution kernels $a_{n-k}^{(n)}$.

\begin{lemma}\label{lem: DCC explicit formula}
For fixed $n\ge2$, the DCC kernels $p_{n-k}^{(n)}$ defined in \eqref{eq: summing orthogonal procedure} satisfy
\beqyy 
p_{n-k}^{(n)}=\frac{1}{a_{0}^{(k)}}-\psi^{(k+1)}_1\sum_{j=k+1}^{n}\frac{1}{a_{0}^{(j)}}
\prod_{\ell=k+2}^{j}\brab{\chi^{(j-\ell)}_{2}\psi_{2}^{(\ell)}-\psi_{1}^{(\ell)}}
\quad\text{for $1\leq k\leq n$,}
\eeqyy
where $\psi^{(m)}_{j}$ and $\chi_{\ell}^{(k)}$ are defined by
Definitions \ref{def: auxiliary sequence Psi} and Definitions \ref{def: auxiliary sequence Chi}, respectively.
\end{lemma}

\begin{proof} Using together Lemma \ref{lem: DOC explicit formula} and \eqref{eq: orthogonal recursive procedure}, we have
\beqyy
\theta_{j-k}^{(j)}=-\frac{1}{a_{0}^{(j)}}\psi^{(k+1)}_1\prod_{\ell=k+2}^{j}\brab{\chi^{(j-\ell)}_{2}\psi_{2}^{(\ell)}-\psi_{1}^{(\ell)}}
\quad\text{for $j\ge k+1$.}
\eeqyy
Thus the definition \eqref{eq: summing orthogonal procedure} yields the desired formula of DCC kernels $p_{n-k}^{(n)}$, namely,
\beyy
p_{n-k}^{(n)}&=&\theta_{0}^{(k)}+\sum_{j=k+1}^{n}\theta_{j-k}^{(j)}\nn \\
&=& \frac{1}{a_{0}^{(k)}}-\psi^{(k+1)}_1\sum_{j=k+1}^{n}\frac{1}{a_{0}^{(j)}}
\prod_{\ell=k+2}^{j}\brab{\chi^{(j-\ell)}_{2}\psi_{2}^{(\ell)}-\psi_{1}^{(\ell)}}
\quad\text{for $1\leq k\leq n$.}
\eeyy
The proof is completed.
\end{proof}

\begin{lemma}\label{lem: n decreasing of DCC kernels}
For any $n\ge2$,
\begin{itemize}
\item
if the discrete convolution kernels $a_{n-k}^{(n)}$ satisfy \textbf{C1}-\textbf{C3},
then the DCC kernels $p_{n-k}^{(n)}$ in \eqref{eq: summing orthogonal procedure} satisfy
\beq \label{2eqp2}
p_{0}^{(n)}>0,\quad p_{0}^{(n-1)}> p_{1}^{(n)}
\quad \text{and}\quad p_{j-1}^{(n-1)}\ge p_{j}^{(n)}\quad\text{for \;\; $2\le j\le n-1$};
\eeq
\item
if the discrete convolution kernels $a_{n-k}^{(n)}$ satisfy the condition \textbf{C4},
then the DCC kernels $p_{n-k}^{(n)}$ in \eqref{eq: summing orthogonal procedure} are non-negative, i.e.,
\beq \label{2eqp3}
p_{n-k}^{(n)}=\sum_{j=k}^{n}\theta_{j-k}^{(j)}\ge 0, \quad\text{for \;\; $1\leq k\le n$.}
\eeq
\end{itemize}
\end{lemma}

\begin{proof} Using Lemma \ref{lem: sign of DOC kernels} under the conditions \textbf{C1}-\textbf{C3}, the results in (\ref{2eqp2}) follow from
the relationship \eqref{eq: orthogonal-complementary relation}:
\[
p_{1}^{(n)}=p_{0}^{(n-1)}+\theta_{1}^{(n)}\le p_{0}^{(n-1)}\quad \text{and}\quad
p_{j}^{(n)}=p_{j-1}^{(n-1)}+\theta_{j}^{(n)}\le p_{j-1}^{(n-1)}\quad\text{for $2\leq j\le n-1.$}
\]
To prove (\ref{2eqp3}), we first take $k=\ell$ and $k=\ell+1$ in the complementary identity \eqref{eq: new complementary identity}
and find the following recursive procedure
\bey\label{eq: complementary recursive procedure}
&& p_{0}^{(n)}=\frac{1}{a^{(n)}_{0}}, \nn \\
&& p_{n-\ell}^{(n)}=\frac{1}{a^{(\ell)}_{0}}\sum_{j=\ell+1}^{n}\brab{a^{(j)}_{j-\ell-1}-a^{(j)}_{j-\ell}}p_{n-j}^{(n)}
\quad\text{for $1\leq \ell\le n-1$.}
\eey
Moreover, the condition \textbf{C4} implies
$a^{(j)}_{j-\ell-1}-a^{(j)}_{j-\ell}\ge0\quad\text{for $j\ge \ell+1$}.$
Thus a simple induction for (\ref{eq: complementary recursive procedure}) yields the desired result (\ref{2eqp3}) .
\end{proof}

We close this section by pointing out that the idea of using the DCC kernels $p_{n-j}^{(n)}$
was first proposed in \cite{LiaoLiZhang:2018,LiaoMcLeanZhang:2019}
for studying the fractional discrete  Gr\"{o}nwall inequalities. On the other hand,
the present work seems the first effort in establishing
the connection between the DCC and DOC kernels for a general class of discrete convolution kernels.
The interplay results for the DOC and DCC kernels can be illustrated by Figure 1.

\begin{figure}[!hbt]
  \centering
  \includegraphics[width=4in]{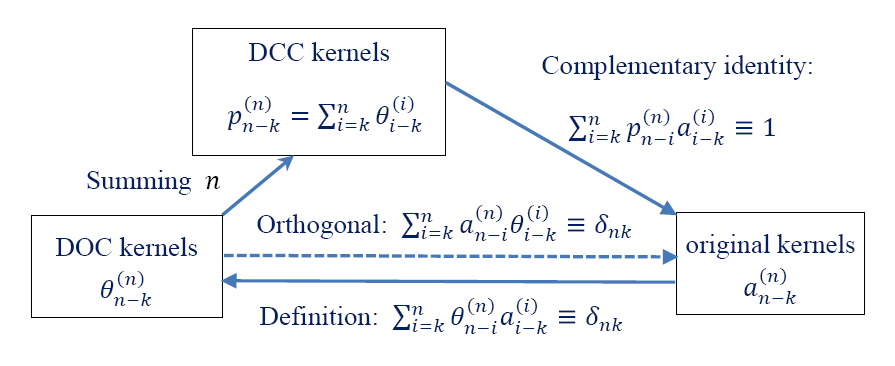}
  \label{fig: DOC and DCC relation}
  \caption{\tt The relationship diagram of DOC and DCC kernels.}
\end{figure}

\section{Proof of Theorem \ref{thm: positive definite quadratic form}} \label{sect3}

 With the preparations of Section 2, we are now ready to prove Theorem \ref{thm: positive definite quadratic form}.

\begin{proof}[Proof of Theorem \ref{thm: positive definite quadratic form}]
It follows from (\ref{2eqc17}) and the definition \eqref{eq: summing orthogonal procedure} that
\beq \label{3eq1}
2\sum_{k=1}^nV_k\sum_{j=1}^k\theta_{k-j}^{(k)}V_j
\ge\sum_{k=1}^n\braB{\sum_{j=1}^k\theta_{k-j}^{(k)}+p_{n-k}^{(n)}}V_k^2\quad \text{for $n\ge1$.}
\eeq
Using (\ref{2eqc15}) and (\ref{2eqp3}) yields the following positivity result:
\beqyy
\sum_{j=1}^k\theta_{k-j}^{(k)}+p_{n-k}^{(n)} >0.
\eeqyy
The above result, together with (\ref{3eq1}),  confirms that the DOC kernels $\theta_{n-k}^{(n)}$ are positive definite.
Then the desired result of Theorem \ref{thm: positive definite quadratic form} follows by using (\ref{2eqc0}).
\end{proof}

Following the proof of Theorem \ref{thm: positive definite quadratic form}, one can easily verify the following result.

\begin{corollary}\label{corol: positive definite quadratic form}
For fixed $n\ge2$, if the discrete convolution kernels $\big\{a^{(n)}_{n-k}\big\}_{k=1}^n$ satisfy
\beqyy
a^{(n)}_{j}\ge0,\quad a^{(n-1)}_{j-1}\ge a^{(n)}_{j},\quad a^{(n-1)}_{j-1}a^{(n)}_{j+1}\ge a^{(n-1)}_{j}a^{(n)}_{j},
\quad a^{(n)}_{j-1}\ge a^{(n)}_{j},
\eeqyy
then the convolution kernels $a^{(n)}_{n-k}$ are positive semi-definite in the sense that
\beqyy
\sum_{k=1}^nw_k\sum_{j=1}^ka^{(k)}_{k-j}w_j\ge0\qquad\text{for any sequence $\{w_{k}\}_{k=1}^n$.}
\eeqyy
\end{corollary}

\begin{proof}
We have $a^{(n)}_{0}>0$ unless the trial case with $a^{(n)}_{j}=0$ for $0\le j\le n-1$.
If the discrete convolution kernels $a^{(n)}_{j}=0$ for $j\ge1$, then the DOC kernels $\theta_{n-k}^{(n)}$ are positive semi-definite
 because $\psi^{(m)}_{\ell}=0$ for any $\ell\ge1$ such that $\theta^{(n)}_{0}>0$ and $\theta^{(n)}_{j}=0$ for $j\ge1$.
 Then using (\ref{2eqc0}) confirms the claimed result.
It remains to consider the case of $a^{(n)}_{0}\ge a^{(n)}_{1}>0$. In this case, similar to the proof
of  the proof for (\ref{2eqc15})-(\ref{2eqc16}) gives
\beyy
&& \theta_{0}^{(n)}>0,\quad  \theta_{j}^{(n)}\le0\;\quad \text{for $1\le j\le n-1$,}\\
&& \sum_{j=1}^n\theta_{n-j}^{(n)}\ge 0. 
\eeyy
Moreover, similar to the proof of (\ref{2eqp3}) gives
\beqyy 
p_{n-k}^{(n)}=\sum_{j=k}^{n}\theta_{j-k}^{(j)}\ge 0, \quad\text{for \;\; $1\leq k\le n$.}
\eeqyy
 Hence, using (\ref{2eqc17}) yields the positive semi-definiteness of
the DOC kernels $\theta_{n-k}^{(n)}$, and then the desired result follows by using Lemma \ref{lem: Mutual orthogonality}
by using (\ref{2eqc0}).
\end{proof}

Corollary \ref{corol: positive definite quadratic form}
presents another class of sufficient conditions for the positive semi-definiteness of discrete kernels $a^{(n)}_{n-k}$.
Notice that if the discrete kernels are independent of the (time-level) index $n$, that is, $a^{(n)}_{j}=a_{j}$,
then the conditions reduce to
\beqyy
a_{j}\ge0,\quad a_{j-1}\ge a_{j},\quad a_{j-1}a_{j+1}\ge a_{j}^2,
\eeqyy
which are slightly stronger than \eqref{cond: LopezMarcos 1990} proposed in \cite{LopezMarcos:1990} since we have
\beqyy
2a_{j}\le 2\sqrt{a_{j-1}a_{j+1}} \le a_{j-1}+a_{j+1}.
\eeqyy

\section{Some applications of Theorem \ref{thm: positive definite quadratic form}} \label{sect4}
\setcounter{equation}{0}

 In this section, we will apply Theorem \ref{thm: positive definite quadratic form} to study
stability of some variable time-stepping schemes.
Recall the definition of the Riemann--Liouville integral operator with
order~$\gamma>0$, (see, e.g.,  \cite{Podlubny:1999}):
\begin{equation}\label{def: Riemann-Liouville integral}
(\mathcal{I}^\gamma v)(t)\defeq\int_0^t\omega_\gamma(t-s)v(s)\zd{s}
	\quad\text{for $t>0$,}
	\quad\text{where $\omega_\gamma(t)\defeq\frac{t^{\gamma-1}}{\Gamma(\gamma)}.$}
\end{equation}
On the other hand, the Caputo fractional derivative of order $\alpha$ is defined by
\begin{equation}\label{def: Caputo derivative}
(\fd{\alpha}v)(t)\defeq(\mathcal{I}^{1-\alpha}v')(t)
	=\int_0^t\omega_{1-\alpha}(t-s)v'(s)\zd{s}
	\quad\text{for $0<\alpha<1$.}
\end{equation}
For a given $T>0$ and the time interval $[0,T]$, we consider a nonuniform time gird
$$0=t_{0}<t_{1}<\cdots<t_{k-1}<t_{k}<\cdots<t_{N}=T$$
with variable step-size $\tau_{k}:=t_{k}-t_{k-1}$ for $1\leq{k}\leq{N}$. The maximum time-step size is defined by $\tau:=\max_{1\leq{k}\leq{N}}\tau_{k}$. Also, we define the local time-step ratio as $r_k:=\tau_k/\tau_{k-1}$ for $k\ge2$.
For a mesh function $v^{k}=v(t_k)$, we set $v^{k-\theta}:=(1-\theta)v^{k}+\theta v^{k-1}$, where $\theta\in[0,1)$ is a weighted parameter.
Also, we define the following difference operators
\beqyy
\diff v^{k}:=v^{k}-v^{k-1} \quad\text{and}\quad\partial_{\tau}v^{k-\frac12}:=\diff v^{k}/\tau_k.
\eeqyy

\subsection{L1 scheme for Caputo fractional derivative}

Our first application is the variable-step L1 scheme for the Caputo fractional derivative \eqref{def: Caputo derivative}.
The L1 formula~(see, e.g. \cite{LiaoLiZhang:2018})
uses $v'(s)\approx \diff v^k/\tau_k$ in the subinterval $(t_{k-1},t_k)$ to obtain
\begin{align*}
(\dfd{\alpha}v)^n:=\sum_{k=1}^nc^{(n,\alpha)}_{n-k}\diff v^k\quad\text{for $n\ge1$,}
\end{align*}
where the associated discrete L1 kernels read
\begin{align}\label{def: L1 kernels}
c_{n-k}^{(n,\alpha)}:=\frac1{\tau_k}\int_{t_{k-1}}^{t_k}\omega_{1-\alpha}(t_n-s)\zd{s}\quad\text{for $1\le k\le n$,}
\end{align}
By using the integral and differential mean value theorems, we have the following result.

\begin{proposition}\label{prop: positive definite L1 kernels}
The discrete L1 kernels $c_{n-k}^{(n,\alpha)}$ in \eqref{def: L1 kernels} satisfy
\begin{align*}
c^{(n,\alpha)}_{j}>0,\quad c^{(n,\alpha)}_{j-1}>c^{(n,\alpha)}_{j},\quad c^{(n-1,\alpha)}_{j-1}>c^{(n,\alpha)}_{j},\quad
c^{(n-1,\alpha)}_{j-1}c^{(n,\alpha)}_{j+1}>c^{(n-1,\alpha)}_{j}c^{(n,\alpha)}_{j}.
\end{align*}
Then the discrete L1 kernels $c^{(n,\alpha)}_{n-k}$ are positive definite according to Theorem \ref{thm: positive definite quadratic form}.
\end{proposition}

\begin{proof} For any fixed time-level index $n\ge2$, the positivity and the decreasing property of $c^{(n,\alpha)}_{j}$ follow
from the following fact which is due to the integral mean value theorem:
$$c^{(n,\alpha)}_{n-k}=\omega_{1-\alpha}(t_n-s_{nk})\quad\text{ for some $s_{nk}\in[t_{k-1},t_k]$}.$$
Consider an auxiliary sequence
\beqyy
\psi_{n-k}^{(n,\alpha)}:=\frac{c_{n-k}^{(n,\alpha)}}{c_{n-1-k}^{(n-1,\alpha)}}
=\frac{c_{n,k}(1)}{c_{n-1,k}(1)}
\quad\text{for $1\le k\le n-1$,}
\eeqyy
where
\begin{align*}
c_{n,k}(\mu):=\frac{1}{\tau_k}\int_{t_{k-1}}^{t_{k-1}+\mu\tau_k}\omega_{1-\alpha}(t_n-s)\zd{s}\quad\text{for  $1\leq k\leq n$.}
\end{align*}
Differentiating it with $\mu$ gives  $c'_{n,k}=\omega_{1-\alpha}(t_n-t_{k-1}-\mu\tau_k)$.
It follows from the Cauchy mean-value theorem that there exists some $\xi_{1k}\in(0,1)$ such that
\bey
\psi_{n-k}^{(n,\alpha)}&=& \frac{c_{n,k}(1)-c_{n,k}(0)}{c_{n-1,k}(1)-c_{n-1,k}(0)}
=\frac{c_{n,k}'(\xi_{1k})}{c_{n-1,k}'(\xi_{1k})}\nn \\
&=& \frac{\omega_{1-\alpha}(t_n-t_{k-1}-\xi_{1k}\tau_k)}{\omega_{1-\alpha}(t_{n-1}-t_{k-1}-\xi_{1k}\tau_k)}\nn \\
&=& \braB{\frac{t_{n-1}-t_{k-1}-\xi_{1k}\tau_k}{t_{n}-t_{k-1}-\xi_{1k}\tau_k}}^{\alpha}
\qquad\text{for $1\leq k\leq n-1$}.
\eey
Note that the function $y=(A-x)/(B-x)$ is decreasing with respect to $x>0$ for $A<B$. Consequently,
\begin{align*}
&\braB{\frac{t_{n-1}-t_{k}}{t_{n}-t_{k}}}^{\alpha}<\psi_{n-k}^{(n,\alpha)}
<\braB{\frac{t_{n-1}-t_{k-1}}{t_{n}-t_{k-1}}}^{\alpha}\quad\text{for $1\leq k\leq n-1$,}
\end{align*}
which yields
\begin{align*}
0<\psi_{1}^{(n,\alpha)}<\psi_{2}^{(n,\alpha)}<\cdots<\psi_{n-1}^{(n,\alpha)}<(t_{n-1}/t_n)^{\alpha}<1\quad\text{for any fixed $n\ge 2$.}
\end{align*}
This leads to the last two inequalities, and the proof is completed.
\end{proof}

\begin{remark}\label{remark: L1+ formula Caputo}
{\rm Recently, a second order L1$^{+}$ formula was proposed in \cite{JiLiaoGongZhang:2019} for the Caputo derivatives:
\beqyy
(\partial_{\tau}^{\alpha}v)^{n-\frac{1}{2}}
:=\frac{1}{\tau_{n}}\int_{t_{n-1}}^{t_{n}}
\int_{0}^{t}\omega_{1-\alpha}(t-s)(\Pi_{1}v)'(s)\zd{s}\zd{t}
=\sum_{k=1}^{n}\bar{c}_{n-k}^{(n)}\triangledown_{\tau}v^{k}\quad \text{for $n\ge1$,}
\eeqyy
where the discrete L1$^{+}$ kernels $\bar{c}_{n-k}^{(n)}$ are defined by
\beq\label{New-L1-Coeff}
\bar{c}_{n-k}^{(n)}
:=\frac{1}{\tau_{n}\tau_{k}}\int_{t_{n-1}}^{t_{n}}
\int_{t_{k-1}}^{\min\{t,t_{k}\}}\omega_{1-\alpha}(t-s)\zd{s}\zd{t}\quad \text{for $1\leq{k}\leq{n}$.}
\eeq
Note that the  positive semi-definiteness of the continuous kernel $\omega_{1-\alpha}(t-s)$
directly leads to the positive semi-definite  property of the discrete L1$^{+}$ kernels $\bar{c}_{n-k}^{(n)}$ \cite[Lemma 3.1]{JiLiaoGongZhang:2019}. However, it can be verified that the condition $\bar{c}_{0}^{(n)}\ge\bar{c}_{1}^{(n)}$ is may fail for some index $n$. In other words,  Theorem \ref{thm: positive definite quadratic form} only provides some easy-to-check sufficient
conditions. There is room for improvement, and it will be ideal to find sufficient and necessary conditions.
}
\end{remark}

\subsection{Energy stability for the time fractional Allen-Cahn equation}

According to Proposition \ref{prop: positive definite L1 kernels}, on general nonuniform time meshes, Theorem 1.1 implies
\begin{align}\label{ieq: positive definite L1 kernels}
\sum_{k=1}^n\diff v^k(\dfd{\alpha}v)^k=\sum_{k=1}^n\diff v^k\sum_{j=1}^kc^{(k,\alpha)}_{k-j}\diff v^j>0\quad\text{if $\diff v^k\not\equiv0$.}
\end{align}
This result can be applied to study the energy stability of variable time stepping schemes for
the time-fractional Allen-Cahn equation \cite{TangYuZhou:2019sisc}:
\begin{align*}
\partial_{t}^{\alpha}u=\varepsilon^{2}\Delta{u}-F'(u)\quad\text{for $x\in\Omega$ and $0<t \le T$,}
\end{align*}
where $\Omega$ is a bounded spatial domain, $\Delta$ is the Laplacian operator and $F(u)=\frac14(1-u^2)^2$ is the well known double well potential.  It is known that for the time-fractional Allen-Cahn equation, there holds the following energy stability  (see, e.g. \cite{Tang20,TangYuZhou:2019sisc}):
\begin{equation*}
  E[u(t)] \leq  E[u(0)] \quad\text{with}\quad  E[u] := \int_{\Omega} \Bigl( \varepsilon^2 |\nabla u|^2 +  F(u)\Bigr)\zd{x}.
\end{equation*}
It is thus desired to require that this stability holds at the discrete level. The following first-order explicit-implicit variable stabilization scheme was proposed in \cite{JiLiaoZhang:2019}:
\begin{align}\label{eq: stabilized scheme Allen-Cahn}
\brab{\partial_{\tau}^{\alpha}u}^{n}
&=\varepsilon^{2}D_h u^{n}
-F'(u^{n-1})-S(u^{n}-u^{n-1})\quad\text{for $n\geq{1}$,}
\end{align}
where $D_{h}$ denotes the discrete matrix of central difference approximation for Laplace operator $\Delta$ with periodic boundary conditions,
and $S(u^{n}-u^{n-1})$ is a stabilization term with $S>0$ being a constant.

It is shown in \cite[Theorem 2.2]{JiLiaoZhang:2019} that scheme \eqref{eq: stabilized scheme Allen-Cahn} preserves
the discrete maximum principle, but the the discrete energy stability proof is not provided. This is partially due to the complexity introduced by the non-uniform time steps. Using the positive definite property \eqref{ieq: positive definite L1 kernels},
one can now easily follow the proof of \cite[Lemma 2.6]{JiLiaoZhang:2019} to obtain the following discrete energy stability result.

\begin{proposition}\label{prop: Allen-cahn Energy-Law}
If the stabilized parameter  $S \ge 2$, then the scheme \eqref{eq: stabilized scheme Allen-Cahn}
preserves the discrete maximum principle, together with the following discrete energy stability
\begin{align*}
E_h^n\le  E_h^0\quad \text{for $n\ge1$},
\end{align*}
where $E_h^n$ is the discrete energy, i.e.,
\begin{align*}
E_h^n=-\frac{\varepsilon^2}{2}(u^n)^TD_{h}u^n
+\sum_{{x}_{h}\in\Omega_{h}}F(u_h^n).
\end{align*}
\end{proposition}
We remark that similar stability results on uniform time meshes has been presented in \cite{JiLiaoZhang:2019} and \cite{TangYuZhou:2019sisc}, while the above stability result on non-uniform grids seems to be new.

\subsection{Numerical Riemann--Liouville integrals}

As an example, we consider the following fractional wave equation
\cite{McLeanMustapha:2007,Podlubny:1999} that involves
the fractional Riemann--Liouville integral:
\begin{equation}\label{eq: fractional wave}
\partial_tu=\mathcal{I}^\gamma \Delta u+f(x,t), \quad x\in\Omega, \quad 0<t \le T.
\end{equation}
It is known that the exact solution admits a weak singularity as $\partial_{tt}u \sim \mathcal{O}(t^{\gamma-1})$ near $t=0$.
One of the effective ways to handle the initial singularity is to use nonuniform time grids, such as the graded meshes \cite{McLeanMustapha:2007,Mustapha:2011,LiaoLiZhang:2018}, to concentrate the time grids near $t=0.$

Similar to \cite{LiXu:2013}, let us consider the backward Euler time-stepping scheme  for \eqref{eq: fractional wave}
\begin{equation}\label{eq: backward scheme}
\partial_{\tau}u^{n-\frac12}=\brab{\mathcal{I}_{\tau}^{\gamma}\Delta u}^{n}+f(x,t_n), \quad n\ge1,
\end{equation}
where $\brab{\mathcal{I}_{\tau}^{\gamma}v}^{n}$ represents the mid-point rule of fractional integral $\mathcal{I}^\gamma$
\begin{align}\label{def: right rectangular RL integral}
\brab{\mathcal{I}_{\tau}^{\gamma}v}^{n}:=\sum_{k=1}^nc^{(n,1-\gamma)}_{n-k}\tau_kv^{k-\frac{1}{2}}, \quad n\ge1.
\end{align}
By Proposition \ref{prop: positive definite L1 kernels}, we see that
\begin{align}\label{ieq: positive definite RL kernels}
\sum_{k=1}^n\tau_kv^{k-\frac{1}{2}}(\mathcal{I}_{\tau}^{\gamma}v)^k
=\sum_{k=1}^n\tau_k v^{k-\frac{1}{2}}\sum_{j=1}^kc^{(k,1-\gamma)}_{k-j} \tau_jv^{j-\frac{1}{2}}>0\quad\text{if $v^k\not\equiv0$.}
\end{align}
By applying property \eqref{ieq: positive definite RL kernels} to (\ref{eq: backward scheme}), one easily gets the following stability estimate
\begin{equation*}
\mynormb{u^{n}}_{L^2}\le \mynormb{u^{0}}_{L^2}+\sum_{k=1}^n\tau_k\mynormb{f(t_k)}_{L^2}.
\end{equation*}
Note that such an $L^2$-norm estimate holds for a general class of nonuniform grids in time. By comparing with the previous stability analysis  \cite{McLeanMustapha:2007,Mustapha:2011,PaniFairweather:2002}, which employed an discrete analogue to the positive semi-definiteness of the continuous kernel, the above analysis provides a more straightforward tool, i.e., by using only the properties of the discrete convolution kernels.

\subsection{Numerical methods for weakly singular Volterra equations}

We next present an application to Volterra integral equations of the form
\begin{equation}\label{eq: Volterra integral}
\partial_tu=\mathcal{K}_t^{(\beta)}\Delta u+f(x,t),\quad x\in\Omega, \quad 0<t \le T.
\end{equation}
In (\ref{eq: Volterra integral}), the convolution integral $\mathcal{K}_t^{(\beta)}v$ is defined by
\begin{equation}\label{def: general convolution integral}
\brab{\mathcal{K}_t^{(\beta)}v}(t):=\int_0^{t}\kappa_{\beta}(t-s)v(s)\zd{s}, \quad t>0,
\end{equation}
where $\kappa_{\beta}(t-s)>0$ is a general positive kernel. We assume that $\kappa_{\beta}'<0$ and $\kappa_{\beta}''\ge 0$ which implies that the kernel $\kappa_{\beta}(t-s)$ is positive definite \cite{McLeanThomeeWahlbin:1996}.

On a general non-uniform grid in time, applying the midpoint rule \cite{Fairweather:1993,McLeanThomeeWahlbin:1996,PaniFairweather:2002} for the convolution integral \eqref{def: general convolution integral} gives
\begin{align*}
\brab{\mathcal{K}_{\tau}^{(\beta)}v}^{k}:=\sum_{k=1}^n\kappa_{n-k}^{(n,\beta)}\tau_kv^{k-\frac{1}{2}},  \quad n\ge1,
\end{align*}
where the discrete kernel is of the form
\begin{align}\label{def: kappa kernels}
\kappa_{n-k}^{(n,\beta)}:=\frac{1}{\tau_k}\int_{t_{k-1}}^{t_{k}}\kappa_{\beta}(t_n-s)\zd{s}, \quad 1\le k\le n.
\end{align}
Note that $y=\kappa_\beta(t_n-x)/\kappa_\beta(t_{n-1}-x)$ is decreasing with respect to $t_{n-1}<x< t_{n}$ . Then one can follow
similar argument as in Proposition \ref{prop: positive definite L1 kernels} to verify the following result.

\begin{proposition}\label{prop: positive definite kappa kernels}
If $\kappa_\beta>0$, $\kappa_{\beta}'<0$ and $\kappa_{\beta}''\ge 0$ , then the discrete convolution kernels $\kappa_{n-k}^{(n,\beta)}$ in \eqref{def: kappa kernels} satisfy
\beqyy
\kappa_{j}^{(n,\beta)}>0,\quad \kappa_{j-1}^{(n,\beta)}>\kappa_{j}^{(n,\beta)},\quad \kappa_{j-1}^{(n-1,\beta)}>\kappa_{j}^{(n,\beta)},\quad
\kappa_{j-1}^{(n-1,\beta)}\kappa_{j+1}^{(n,\beta)}\ge\kappa_{j}^{(n-1,\beta)}\kappa_{j}^{(n,\beta)}.
\eeqyy
Consequently, the discrete kernels $\kappa_{n-k}^{(n,\beta)}$ are positive definite according to Theorem \ref{thm: positive definite quadratic form}.
\end{proposition}

Now, consider the backward Euler time-stepping scheme for solving \eqref{eq: Volterra integral}:
\beq \label{4eq10}
\partial_{\tau}u^{n-\frac12}=\brab{\mathcal{K}_{\tau}^{(\beta)}\Delta u}^{n}+f(x,t_n), \quad n\ge1.
\eeq
By Proposition \ref{prop: positive definite kappa kernels}, we see that
\begin{align}\label{ieq: positive definite Volterra kernels}
\sum_{k=1}^n\tau_kv^{k-\frac{1}{2}}(\mathcal{K}_{\tau}^{(\beta)}v)^k
=\sum_{k=1}^n\tau_k v^{k-\frac{1}{2}}\sum_{j=1}^k\kappa^{(k,\beta)}_{k-j} \tau_jv^{j-\frac{1}{2}}>0\quad\text{if $v^k\not\equiv0$.}
\end{align}
By applying \eqref{ieq: positive definite Volterra kernels} to the backward Euler scheme (\ref{4eq10}), it is easy to obtain the following $L^2$-norm estimate:
\beqyy
\mynormb{u^{n}}_{L^2}\le \mynormb{u^{0}}_{L^2}+\sum_{k=1}^n\tau_k\mynormb{f(t_k)}_{L^2}, \quad n\ge1.
\eeqyy
We emphasize again that this estimate is valid for a general class of non-uniform grid in time, including the well-studied graded meshes.

\begin{remark}\label{remark: CN formula general}
{\rm
For general convolution integrals \eqref{def: general convolution integral}
with smooth or weakly singular kernels, a second-order Crank-Nicolson formula was proposed by McLean et al. \cite{McLeanThomeeWahlbin:1996,McLeanMustapha:2007}:
\beq\label{New-L1-Formula}
(\mathcal{K}_{\tau}^{(\beta)}v)^{n-\frac{1}{2}}
:=\frac{1}{\tau_{n}}\int_{t_{n-1}}^{t_{n}}
\int_{0}^{t}\kappa_{\beta}(t-s)v^{k-\frac12}\zd{s}\zd{t}
=\sum_{k=1}^{n}\bar{\kappa}_{n-k}^{(n)}\tau_kv^{k-\frac12}\quad \text{for $n\ge1$,}
\eeq
where the associated discrete convolution kernels  $\bar{\kappa}_{n-k}^{(n)}$ are defined by
\beq\label{New-L1-CoeffX}
\bar{\kappa}_{n-k}^{(n)}
:=\frac{1}{\tau_{n}\tau_{k}}\int_{t_{n-1}}^{t_{n}}
\int_{t_{k-1}}^{\min\{t,t_{k}\}}\kappa_{\beta}(t-s)\zd{s}\zd{t}\quad \text{for $1\leq{k}\leq{n}$.}
\eeq
Again, the positive semi-definiteness of the continuous kernel
directly implies the positive semi-definite property of the discrete kernels $\bar{\kappa}_{n-k}^{(n)}.$
Similar as in Remark \ref{remark: L1+ formula Caputo}, the condition \textbf{C4} is not always fulfilled here because of the violation of $\bar{\kappa}_{0}^{(n)}\ge \bar{\kappa}_{1}^{(n)}$.
Nevertheless, it is expected that our technique in proving Theorem \ref{thm: positive definite quadratic form}
is still useful for handling this situation, and this will be investigated in a future work.
}
\end{remark}


\section{Concluding remarks}
\setcounter{equation}{0}

The positive definiteness of real quadratic forms with convolution structure plays an important role in analyzing variable time-stepping schemes for time-fractional  differential equations and Volterra-type integral equations. In contrast to the classical treatments which use a discrete analogue to the positive semi-definiteness of continuous kernels,
we present in this work a unified criteria obtained by using  the DOC and DCC tools. More precisely, for convolution
kernels $a^{(n)}_{n-k}$ we show for the first time some  sufficient algebraic conditions that implies the positive
definiteness of the associated real quadratic form. We emphasize the simplicity of the imposed conditions which
are relevant to  positivity (\textbf{C1}), monotonicity (\textbf{C2} and \textbf{C4}), and convexity (\textbf{C3}).
The usefulness and easy-to-check style are demonstrated in Section \ref{sect4} where several first-order schemes are examined.

As noticed in Remarks \ref{remark: L1+ formula Caputo}-\ref{remark: CN formula general}, the condition \textbf{C4}
is not sharp for treating second-order discretization schemes.
Nevertheless, similar situations occur even for the case of uniform grids. For example,
the condition \eqref{cond: LopezMarcos 1990} of \cite{LopezMarcos:1990} works well for first-order schemes but
failed for second-order approximations, see, e.g., \cite{Tang:1993}.
On the other hand, it is expected that the condition \textbf{C4} may be weakened so that our conditions are useful for high order
time discretizations. A careful examination of the proof of Theorem \ref{thm: positive definite quadratic form} reveals that conditions \textbf{C1}-\textbf{C3} together with the following inequality
\begin{align}\label{cond: weaken C4}
p_{n-k}^{(n)}>-\sum_{j=1}^k\theta_{k-j}^{(k)}\quad\text{for $1\le k\le n$}
\end{align}
are sufficient for the positive definiteness claim.
Note that the condition \textbf{C4} ensures a stronger estimate $p_{n-k}^{(n)}\ge0$, see \eqref{2eqp3}.
However, the weaker constraint \eqref{cond: weaken C4} is quite complicated and not as elegant as desired.  Some further
efforts may be made in this direction.


\begin{thebibliography}{99}


\bibitem{ChenXuCaoZhou:2018}
{\sc H.-B. Chen, D. Xu, J. Cao and J. Zhou},
{\em A backward Euler alternating direction implicit difference scheme
for the three-dimensional fractional evolution equation},
Numer Methods Partial Differential Eq., 34 (2018), pp. 938--958.
%
\bibitem{CuestaPalencia:2003}
{\sc E. Cuesta and C. Palencia},
{\em A fractional trapezoidal rule for integro-differential equations of fractional order in Banach spaces},
App. Numer. Math., 45 (2003), pp. 139--159.

\bibitem{Fairweather:1993}
{\sc G. Fairweather},
{\em Spline collocation methods for a class of hyperbolic partial integro-differential equations},
SIAM J. Numer. Anal., 31(2) (1993), pp. 444--460.

\bibitem{GrenanderSzego:1958}
{\sc U. Grenander and G. Szeg\"{o}}, {\em Toeplitz Forms and their Applications},
Second (textually unaltered) edition, Chelsea Publishing Company, New York, 1984
(First edition publicated by University of California Press, Berkeley, CA, 1958).


\bibitem{JiLiaoZhang:2019}
{\sc B. Ji, H.-L. Liao and L. Zhang},
{\em Simple maximum-principle preserving time-stepping methods
for time-fractional Allen-Cahn equation}, Adv. Comput. Math.,
46(2) (2020), doi: 10.1007/s10444-020-09782-2.


\bibitem{JiLiaoGongZhang:2019}
{\sc B. Ji, H.-L. Liao, Y. Gong and L. Zhang},
{\em Adaptive second-order Crank-Nicolson time-stepping schemes for
 time fractional molecular beam epitaxial growth models},
SIAM J. Sci. Comput.,  2020, 42(3) : B738-B760.



\bibitem{LiXu:2013}
{\sc L. Li and D. Xu},
{\em Alternating direction implicit-Euler method for the two-dimensional fractional evolution equation},
J. Comput. Phys., 236 (2013), pp. 157--168.

\bibitem{LiaoLiZhang:2018}
{\sc H.-L. Liao, D. Li and J. Zhang},
{\em Sharp error estimate of nonuniform L1 formula for linear
reaction-subdiffusion equations},
SIAM J. Numer. Anal., 56(2) (2018), pp. 1112-1133.


\bibitem{LiaoMcLeanZhang:2019}
{\sc H.-L. Liao, W. McLean and J. Zhang},
{\em A discrete {Gr\"{o}nwall} inequality
with application to numerical schemes for subdiffusion problems},
SIAM J. Numer. Anal., 57(1)  (2019), pp. 218-237.




\bibitem{LiaoTangZhou:2019}
{\sc H.-L. Liao, T. Tang and T. Zhou},
{\em A second-order and nonuniform time-stepping maximum-principle preserving scheme for time-fractional
Allen-Cahn equations}, J. Comput. Phys., (414) 2020, 109473.


\bibitem{LiaoZhang:2019}
{\sc H.-L. Liao and Z. Zhang},
{\em Analysis of adaptive BDF2 scheme for diffusion equations},
Math. Comput., 2020, DOI: 10.1090/mcom/3585.

%

\bibitem{LinXu:2007}
{\sc Y.~Lin and C.~Xu},
{\em Finite difference/spectral approximations for the time-fractional
  diffusion equation}, J. Comput. Phys., 225(2) (2007), pp. 1533-1552.

\bibitem{LopezMarcos:1990}
{\sc J.C. L\'{o}pez-Marcos},
{\em A difference scheme for a nonlinear partial integr-odifferential  equation},
SIAM J. Numer. Anal., 27 (1990), pp. 20--31.




\bibitem{Lubich_I}
{\sc C. Lubich,}
{\em Convolution quadrature and discretized operational calculus I.}
Numer. Math. 52(1988), 129-145.


\bibitem{Lubich_II}
{\sc C. Lubich,}
{\em Convolution quadrature and discretized operational calculus II.}
Numer. Math. 52(1988), 413-425.


\bibitem{LubichSloanThomee:1996}
{\sc C. Lubich, I.H. Sloan, and V. Thom\'{e}e},
{\em Nonsmooth data error estimates for approximations of an evolution equation
with a positive-type memory term}, Math. Comp., 65 (1996), pp. 1--17.


\bibitem{McLeanThomeeWahlbin:1996}
{\sc W. McLean, V. Thom\'{e}e, L. B. Wahlbin},
{\em Discretization with variable time steps of an evolution
equation with a positive-type memory term},
J. Comput. Appl. Math., 69 (1996), pp. 49-69.

\bibitem{McLeanMustapha:2007}
{\sc W.~McLean and K.~Mustapha},
{\em A second-order accurate numerical method
  for a fractional wave equation},
  Numer. Math., 105 (2007), pp.~481--510.

\bibitem{Mustapha:2011}
{\sc K. Mustapha},
\emph{An implicit finite difference time-stepping method for
a subdiffusion equation with spatial discretization by finite elements},
IMA J. Numer. Anal., 31 (2011), pp.~719-739.


\bibitem{PaniFairweather:2002}
{\sc A. K. Pani and G. Fairweather},
{\em An $H^1$-Galerkin mixed finite element method for an evolution equation with a positive-type memory term},
SIAM J. Numer. Anal., 40(4) (2002), pp.~1475--1490.



\bibitem{Podlubny:1999}
{\sc I. Podlubny}, \emph{Fractional differential equations}, Academic Press, New
York, 1999.

\bibitem{SunWu:2006}
{\sc Z. Sun and X. Wu},
{\em A fully discrete difference scheme for a diffusion-wave system},
App. Numer. Math., 56 (2006), pp. 193--209.

\bibitem{Tang:1993}
{\sc T. Tang},
{\em A finite difference scheme for partial integro-differential equations with a weakly singular kernel},
App. Numer. Math., 11 (4) (1993), pp. 309--319.

\bibitem{Tang20}
{\sc T. Tang},
{\em Revisit of semi-implicit schemes for phase-field equations}, To appear in Analysis in Theory and Applications, 2020.

\bibitem{TangYuZhou:2019sisc}
T.~Tang, H.~Yu, and T.~Zhou.
{\em On energy dissipation theory and numerical stability for
  time-fractional phase field equations},
SIAM J. Sci. Comput., 41 (2019), pp. A3757--A3778.


%
%
%
%

%
%
%
%
%

%
%
%
%
%
%
%
%
%
%
%
%


%
%
%
%


%
%
%
%
%
%
%
%
%

%


%


\end{thebibliography}
\end{document}